\documentclass[12pt,a4]{article}
\usepackage{amsmath}
\usepackage{amscd}
\usepackage{amstext}
\usepackage{amsfonts}
\usepackage{euscript}

\def\scr{\EuScript}

\newcommand{\CC}{{\mathbb C}}
\newcommand{\hol}{{\scr O}}
\newcommand{\HOL}{{\scr O}_{\CC^n}}
\newcommand{\DX}{{\scr D}_X}
\newcommand{\OX}{{\scr O}_X}

\DeclareMathOperator{\der}{\rm Der}

\DeclareMathOperator{\Ann}{\rm Ann}

\DeclareMathOperator{\gr}{\rm Gr}

\DeclareMathOperator{\Gr}{\it Gr}

\DeclareMathOperator{\IM}{\rm Im}

\DeclareMathOperator{\Der}{\it Der}

\DeclareMathOperator{\fHom}{\it Hom}

\DeclareMathOperator{\simm}{Sim}
\newcommand{\DER}{\Der_{\CC}(\HOL)}
\newcommand{\derhol}{\der_{\CC}(\hol)}
\newcommand{\derlogf}{\der (\log f)}
\newcommand{\logD}{\log D}
\newcommand{\derD}{\Der(\logD)}
\newcommand{\derlog}{\Der(\log f)}
\newcommand{\Olog}{\Omega_X^{\bullet}(\logD)}
\newcommand{\GrD}{\Gr_{F^{\bullet}}(\DX)}
\newcommand{\grD}{\gr_{F^{\bullet}}({\scr D})}
\newcommand{\D}{{\scr D}}
\newcommand{\R}{{\cal R}}
\newcommand{\J}{{\scr J}}
\def\O{{\Omega}}
\def\p{{\partial}}
\def\Dx{\p_x}
\def\Dy{\p_y}
\def\Dz{\p_z}
\def\b{{\bullet}}
\newcommand{\depth}{\operatorname{\rm depth}}
\def\bK{{\bf K}}

\newcounter{numero}[section]
\renewcommand{\thenumero}{\thesection .\arabic{numero}}
\newenvironment{corolario}{
\dimen255=\parindent \parindent=0in
\refstepcounter{numero}{\medskip}{\bf Corollary \thenumero.--}\
}{\parindent=\dimen255\par\vspace{1ex}}

\newenvironment{lema}{
\dimen255=\parindent \parindent=0in
\refstepcounter{numero}{\medskip}{\bf Lemma \thenumero.--}\
}{\parindent=\dimen255\par\vspace{1ex}}

\newenvironment{ejemplo}{
\dimen255=\parindent \parindent=0in
\refstepcounter{numero}{\medskip}{\bf Example \thenumero.--}\
}{\parindent=\dimen255\par\vspace{1ex}}

\newenvironment{proposicion}{
\dimen255=\parindent \parindent=0in
\refstepcounter{numero}{\medskip}{\bf Proposition \thenumero.--}\
} {\parindent=\dimen255\par\vspace{1ex}}

\newenvironment{definicion}{
\dimen255=\parindent \parindent=0in
\refstepcounter{numero}{\medskip}{\bf Definition \thenumero.--}\ }
{\parindent=\dimen255\par\vspace{1ex}}

\newenvironment{teorema}{
\dimen255=\parindent \parindent=0in
\refstepcounter{numero}{\medskip}{\bf Theorem \thenumero.--}\ }
{\parindent=\dimen255\par\vspace{1ex}}

\newenvironment{conjetura}{
\dimen255=\parindent \parindent=0in
\refstepcounter{numero}{\medskip}{\bf Conjecture \thenumero.--}\ }
{\parindent=\dimen255\par\vspace{1ex}}

\newenvironment{nota}{
\dimen255=\parindent \parindent=0in
\refstepcounter{numero}{\medskip}{\bf Remark \thenumero.--}\ }
{\parindent=\dimen255\par\vspace{1ex}}

\newenvironment{problema}{
\dimen255=\parindent \parindent=0in
\refstepcounter{numero}{\medskip}{\bf Problem \thenumero.--}\ }
{\parindent=\dimen255\par\vspace{1ex}}

\newcommand\numero{\medskip\refstepcounter{numero}\noindent{\bf
 \thenumero}\hspace{1em}}

\newenvironment{prueba}{{\bf Proof:}\ }{{\hfill C.Q.D.}
\par\vspace{2ex}}
\title{The module ${\scr D}f^s$ for locally quasi-homogeneous free divisors}

\author{Francisco Calder\'{o}n-Moreno and Luis Narv\'{a}ez-Macarro\thanks{
The authors are supported by PB97-0723.}}
\date{}

\begin{document}

\maketitle 

\pagenumbering{arabic} \pagestyle{plain}

 {\small
\begin{center} {\bf Abstract} \end{center}

 We find explicit free resolutions for the $\D$-modules ${\scr D}
f^s$ and ${\scr D}[s] f^s/{\scr D}[s] f^{s+1}$, where $f$ is a
reduced equation of a locally quasi-homogeneous free divisor.
These results are based on the fact that every locally
quasi-homogeneous free divisor is Koszul free, which is also
proved in this paper.}

\section*{Introduction}

In this paper we study the module ${\scr D} f^s$, where $\scr D$
is the ring of germs at $0\in \CC^n$ of linear holomorphic
differential operators and $f$ is a reduced local equation of a
locally quasi-homogeneous free divisor $D\subset (\CC^n,0)$.

The module ${\scr D} f^s$ encodes an enormous amount of geometric
information of the singularity $f=0$, but usually it is hard to
work with in an explicit way. We prove the following results (see
corollary \ref{coro-x} and theorem \ref{mainn}):

\medskip\noindent
{\bf (A)} Let $f=0$ be a reduced local equation of a locally
quasi-homogeneous free divisor of $\CC^n$,  and let
$\{\delta_1,\dots,\delta_{n-1}\}$ be a basis of the module of
vector fields vanishing on $f$. Then
\begin{enumerate}
\item
The $\delta_i$ generate the ideal  $\Ann_{\scr D} f^s$.
\item There exist  explicit Koszul-Spencer type free resolutions
for the modules ${\scr D} f^s$ and ${\scr D}[s] f^s/{\scr D}[s]
f^{s+1}$ built on $\delta_1,\dots,\delta_{n-1}$ and
$f,\delta_1,\dots,\delta_{n-1},$ respectively.
\end{enumerate}

\medskip
Locally quasi-homogeneous free divisors form an important class of
divisors with non isolated singularities: normal crossing
divisors, the union of reflecting hyperplanes of a complex
reflection group, free hyperplane arrangements or the discriminant
of stable mappings in Mather's ``nice dimensions" are examples of
such divisors.

\medskip Let $X$ be a complex analytic manifold. Given a divisor $D\subset X$,
let us write $j:U=X\setminus D\hookrightarrow X$ for the
corresponding open inclusion and $\O^\b(\ast D)$ for the
meromorphic de Rham complex with poles along $D$. In \cite{grot},
Grothendieck proved that the canonical morphism $\O^\b(\ast D) \to
{\bf R} j_\ast (\CC_U)$ is an isomorphism (in the derived
category). This result is usually known as (a version of) {\em
Grothendieck's Comparison Theorem}.

In \cite{ksaito_log}, K. Saito introduced the {\em logarithmic de
Rham complex} associated with $D$, $\Olog$, generalizing the well
known case of normal crossing divisors (cf. \cite{deligne}). In
the same paper, K. Saito also introduced the important notion of
{\em free divisor}.

In \cite{cas_mon_nar_96}, it is proved that  the logarithmic de
Rham complex $\Olog$ computes the cohomology of the complement $U$
if $D$ is a locally quasi-homogeneous free divisor (we say that
$D$ satisfies the {\em logarithmic comparison theorem}). In other
words, the canonical morphism $\Olog \to {\bf R} j_\ast (\CC_U)$
is an isomorphism, or using Grothendieck's result, the inclusion
$\Olog \hookrightarrow \O^\b(\ast D)$ is a quasi-isomorphism. In
fact, in \cite{cal_mon_nar_cas} it is proved that, in the case of
$\dim X = 2$, $D$ is locally quasi-homogeneous if and only if it
satisfies the logarithmic comparison theorem.

Since the derived direct image ${\bf R} j_\ast (\CC_U)$  is a
perverse sheaf (it is the de Rham complex of the  holonomic module
of meromorphic functions with poles along $D$ \cite{mebkhout}, II,
th. 2.2.4), we deduce that $\Olog$ is perverse for every locally
quasi-homogeneous free divisor.

On the other hand,  the first author proved the following results
\cite{cal_99}:\\ Let $D\subset X$ be a Koszul free divisor (see
definition \ref{kfd}) and $\scr I$ the
 left ideal of the ring $\DX$ of differential operators on $X$
 generated by the logarithmic vector fields with respect to $D$. Then\\
1) The left
 $\DX$-module $\DX/{\scr I}$ is holonomic.\\
 2) There is a canonical isomorphism in the derived category
$$\Omega_X^{\bullet}(\logD) \simeq {\bf R} \fHom_{{\scr
D}_X}( \DX/{\scr I},
 \OX  ). $$
 As a consequence, the logarithmic de Rham complex
 associated with a Koszul free divisor is a perverse sheaf.

The proof of {\bf (A)} depends strongly on the following result,
which has been suggested by the above results (see theorem
\ref{lqhiskfd}):

\medskip\noindent
{\bf (B)} Every locally quasi-homogeneous free divisor is Koszul
free.

\medskip
In the first three sections we introduce some material concerning
locally quasi-homogeneous free divisors, Koszul free divisors, the
notion of linear type ideal and the module ${\scr D}f^s$.

In the fourth section we include the proof of {\bf (B)} in our
previous paper \cite{lqhf-kf}.

The fifth section is the main part of this paper and contains the
proof of {\bf (A)} and some related results.

In the sixth section we study some examples and we state some
problems and conjectures.

\medskip
The first part of {\bf (A)} has been proposed (without proof) in
\cite[page 240]{alek_90} in the particular case of discriminants
of versal deformations of simple hypersurface singularities. The
normal crossing divisors case has been treated in \cite{gros_nar}.

\section{Locally quasi-homogeneous and Koszul free divisors}
\numero \label{notations}
 Let $X$ be a $n$-dimensional complex analytic manifold.
 We denote by  $\pi:
 T^{\ast}X \to X$ the cotangent bundle, $\OX$ the sheaf of
 holomorphic functions on $X$,
 $\DX$ the sheaf of linear
 differential operators on $X$ (with holomorphic coefficients),
 $\GrD$ the graded ring associated with the filtration $F^{\bullet}$ by the
 order, $\sigma(P)$ the principal symbol of a differential
 operator $P$ and $\{-,-\}$ the Poisson bracket on ${\scr O}_{T^{\ast}X}$ or $\GrD$.
  We will note ${\scr O}={\scr O}_{X,p}$, ${\scr D}={\scr D}_{X,p}$
  and $\grD={\GrD}_p$
  the respective stalks at $p$, with $p$ a point of $X$. If
  $J\subset \D$ is a left ideal, we denote by $\sigma (J)$ the
  corresponding graded ideal of $\grD$.
 Given a divisor $D\subset X$, we denote by $\derD$ the $\OX$-module of
 the logarithmic vector fields with respect to $D$
 \cite{ksaito_log}. If $f$ is a local equation of $D$ at $p$, we
 denote by $\derlogf$ the stalk at $p$ of $\derD$, whose elements
 are germs at $p$ of vector fields $\delta$ such that $\delta (f)
 \in (f)$.

\begin{definicion} \label{eh}
A divisor $D$  is Euler-homogeneous at $p\in D$ if there is a
local equation $h$ for $D$ around $p$, and a germ of (logarithmic)
vector field $\delta$ such that $\delta(f)=f$. A such $\delta$ is
called a local Euler vector field for $f$.
\end{definicion}

The set of points where a divisor is Euler-homogeneous is open.

\begin{definicion} (cf. \cite{cas_mon_nar_96}) \label{lqh}
A germ of divisor $(D,p)\subset (X,p)$ is quasi-homogeneous if
there are local coordinates  $x_1,\dots,x_n\in\hol_{X,p}$ with
respect to which $(D,p)$ has a weighted homogeneous defining
equation (with strictly positive weights). A divisor $D$ in a
$n$-dimensional complex manifold $X$ is locally quasi-homogeneous
if the germ $(D,p)$ is quasi-homogeneous for each point $p\in D$.
A germ of divisor $(D,p)\subset (X,p)$ is locally
quasi-homogeneous if the divisor $D$ is locally quasi-homogeneous
in a neighborhood of $p$.
\end{definicion}

Obviously a locally quasi-homogeneous divisor is Euler-homogeneous
at every point.

\begin{definicion} \label{def-1} We say that a reduced germ $f\in \hol_{X,p}$ is
locally quasi-homogeneous if the germ of divisor $(\{f=0\},p)$ is.
\end{definicion}

\begin{nota} A reduced germ $f\in \hol_{X,p}$ is
locally quasi-homogeneous if and only if for every $q\in \{f=0\}$
near $p$ there are local coordinates $z_1,\dots,z_n\in\hol_{X,q}$
and a quasi-homogeneous polynomial $P(t_1,\dots,t_n)$ (with
strictly positive weights) such that $f_q= P(z_1,\dots,z_n)$.
\end{nota}

\begin{definicion} \label{kfd}
(\cite{ksaito_log}, \cite{cal_99}, def. 4.1.1) Let $D\subset X$ be
a divisor. We say that $D$ is free at $p\in X$ if $\derD_p$ is a
free $\scr O$-module (of rank $n$). We say that $D$ is a Koszul
free divisor at $p\in X$ if it is free at $p$ and there exists a
basis $\{\delta_1,\dots,\delta_n\}$ of $\derD_p$ such that the
sequence of symbols $\{\sigma(\delta_1),\dots,\sigma(\delta_n)\}$
is regular in $\grD=\GrD_p$ . If $D$ is a free (resp. Koszul free)
divisor at each point of $X$, we simply say that it is free (resp.
Koszul free). We say that a  reduced germ $f\in \hol_{X,p}$ is
free if the divisor $f^{-1}(0)$ is free at $p$.

\end{definicion}

Let's remark that a divisor $D$ is automatically Koszul free at
every $p\in X\setminus D$.

\begin{nota}
The ideal $I_{D,p}=\grD\derD_p$ is generated by the elements of
any basis of $\derD_p$. As $D$ is Koszul free at $p$ if and only
if $\depth(I_{D,p},\grD)=n$ (cf. \cite{mat}, cor.~16.8),
 it is
clear that the definition of Koszul free divisor does not depend
on the election of a particular basis.
 By the coherence of $\GrD$,
 if a divisor is Koszul free at a point,
then it is Koszul free near that point.
\end{nota}
We have not found a reference for the following well known
proposition (see \cite{mat}, th. 17.4 for the local case).

\begin{proposicion} \label{equiv} Let $\CC\{x\}$ be the ring of
convergent power series in the variables $x = x_1,\dots,x_n$ and
let $G$ be the graded ring of polynomials in the variables
$\xi_1,\dots,\xi_t$ with coefficients in $\CC\{x\}$. A sequence
$\sigma_1, \dots, \sigma_s$ of homogeneous polynomials in $G$ is
regular if and only if the set of zeros $V(I)$ of the ideal $I$
generated by $\sigma_1,\dots,\sigma_s$ has dimension $n+t-s$ in
$U\times\CC^t$, for some open neighborhood $U$ of 0 (then each
irreducible component has dimension $n+t-s$).
\end{proposicion}
\begin{prueba}
Let $\CC\{x,\xi\}$ be the ring of convergent power series in the
variables $ x_1,\dots,x_n,\xi_1,\dots,\xi_t$. As the $\sigma_i$
are homogeneous in $G$ and the ring $\CC\{x,\xi\}$ is a flat
extension of $G$, the $\sigma_i$ are a regular sequence in $G$ if
and only if they are a regular sequence in $\CC\{x,\xi\}$. But the
last condition is equivalent to the equality ({\it loc.~cit.}): $$
\dim_{(0,0)}(V(I)) =  \dim \left(\CC\{x,\xi\}/I\right)=n+t-s.$$
Finally, using the fact that all the $\sigma_i$ are homogeneous in
the variables $\xi$, the local dimension of  $V(I)$ at $(0,0)$
coincides with its dimension in $U\times\CC^t$ for some
neighborhood $U$ of 0.
\end{prueba}
\begin{corolario} \label{2.9}
Let $D\subset X $ be a free divisor. Let $J$ be the ideal in
${\scr O }_{T^{\ast} X}$ generated by $\pi^{-1}\derD$. Then, $D$
is Koszul free if and only if the set $V(J)$ of zeros of $J$ has
dimension $n$ (in this case, each irreducible component of $V(J)$
has dimension $n$).
\end{corolario}
\begin{proposicion} \label{bases}
 Let $X$ be a complex manifold of dimension $n$ and let $D\subset
X$ be a divisor. Then:
\begin{enumerate}
\item  Let $X^{\prime}=X\times\CC$ and $D^{\prime}=D\times\CC$.
The divisor $D\subset X$ is Koszul free if and only if
$D^{\prime}\subset X^{\prime}$ is Koszul free.
\item
 Let $Y$ be another complex manifold of dimension $r$ and let $E\subset
Y$ be a divisor. Then: a) The divisor $(D\times Y) \cup (X\times
E)$ is free if $D\subset X$ and $E\subset Y$ are free.
\\
b) The divisor $(D\times Y) \cup (X\times E)$ is Koszul free if
$D\subset X$ and $E\subset Y$
 are Koszul free.
\end{enumerate}
\end{proposicion}
\begin{prueba}
\begin{enumerate}
\item It is a consequence of
\cite{cas_mon_nar_96}, lemma 2.2, (iv) and the fact that
$\sigma_1,\dots,\sigma_n$ is a regular sequence in
$\hol_{X,p}[\xi_1,\dots,\xi_n]$ if and only if
$\xi_{n+1},\sigma_1,\dots,\sigma_n$ is a regular sequence in
$\hol_{X^{\prime},(p,t)}[\xi_1,\dots,\xi_n,\xi_{n+1}]$.
\item a) It is an immediate consequence of Saito's criterion (cf. \cite{cas_mon_nar_96}, lemma 2.2, (v)). \\
      b) It is a consequence of a) and Corollary \ref{2.9}.
\end{enumerate}
\end{prueba}
\begin{ejemplo} \label{examples} Examples of Koszul free divisors are: \\
1) Nonsingular divisors.\\ 2) Normal crossing divisors. \\ 3)
Plane curves: If dim$_{\CC}X=2$, we know that every divisor
$D\subset X$ is free \cite{ksaito_log}, cor. 1.7. Let
$\{\delta_1,\delta_2\}$ be a basis of $\derD_x$. Their symbols
$\{\sigma_1,\sigma_2\}$ are obviously linearly independent over
$\hol$, and by  Saito's criterion \cite{ksaito_log}, 1.8, they are
relatively primes in $\grD=\hol[\xi_1,\xi_2]$. So they form a
regular sequence in $\grD$, and $D$ is Koszul free (see
\cite{cal_99}, cor. 4.2.2).\\ 4) Proposition \ref{bases} gives a
way to obtain Koszul free divisors in any dimension.\\ 5) There
are irreducible Koszul free divisors in dimensions greater than 2,
which are not constructed from divisors in lower dimension
\cite{saito}:
  $X=\CC^3$ and $D\equiv\{f=0\}$,
with $$f=2^8z^3-2^7x^2z^2+2^4x^4z+2^43^2xy^2z-2^2x^3y^2-3^3y^4.$$
A basis of $\derlog$ is $\{\delta_1,\delta_2,\delta_3\}$, with $$
\begin{array}{cccccccccc}
\delta_1&=&6y&\partial_x &+& (8z-2x^2)
&\partial_y&-&xy&\partial_z, \\ \delta_2&=&(4x^2-48z)&\partial_x
&+& 12x y &\partial_y&+&(9y^2-16xz)&
\partial_z, \\
\delta_3&=&2x&\partial_x &+&3y&\partial_y& +& 4z&\partial_z,
\end{array}
$$ and the sequence
$\{\sigma(\delta_1),\sigma(\delta_2),\sigma(\delta_3)\}$ is
$\grD$-regular.
\end{ejemplo}

\section{Ideals of linear type}

\begin{definicion} \label{ilt} (cf. \cite{vascon_cmcaag}, \S 7.2) Let $A$ be a commutative
ring, $I\subset A$ an ideal, $\R (I)= \oplus_{i=0}^\infty I^d t^d
\subset A[t]$ the Rees algebra of $I$ and
 $\simm(I)$ the
symmetric algebra of the $A$-module $I$. We say that $I$ is of
{\it linear type} if the canonical (surjective) morphism of graded
$A$-algebras $$ \simm(I) \to \R(I)$$ is an isomorphism.
\end{definicion}

\begin{lema} \label{kerishom}
 Given a commutative ring $A$ and an ideal $I\subset A$ generated by a family of
elements $\{a_i\}_{i\in\Lambda}$, the following properties are
equivalent:
\begin{enumerate}
\item[a)] $I$ is of linear type.
\item[b)] If $\varphi:A[\{X_i\}_{i\in\Lambda}]\to \R(I)$ is the morphism of
graded algebras defined by $\varphi(X_i) = a_it$, then the kernel
of $\varphi$ is generated by homogeneous elements of degree 1.
\end{enumerate}
\end{lema}

\begin{prueba}
We consider the kernel of the surjective morphism of graded
$A$-algebras $ \Phi: A[\{X_i\}_{i\in\Lambda}]\to \simm(I),$
defined by $\Phi(X_{i_1}\cdots X_{i_d}) = a_{i_1}\cdot \ldots
\cdot a_{i_d}.$ Then $\ker(\Phi)= \ker (\varphi)$ if and only if
$I$ is of linear type, $\ker(\Phi)$ is an ideal generated by its
homogeneous elements of degree 1, $\ker(\Phi)_1$, and
$\ker(\Phi)_1 = \ker(\varphi)_1$.
\end{prueba}

The definition and the lemma above sheafify in the obvious way.

\medskip
The following results concern the case where the ideal $I$ is
generated by a regular sequence.

\begin{lema} \label{regular2}
Let $\{a_1,\cdots,a_m\}$ be an $A$-sequence. For $p\leq m$, if
$\alpha a_1^{s_1}\cdots a_{m}^{s_{m}}\in( a_1^{s_1+k_1},\cdots,
a_{p}^{s_p+k_p})$, then $\alpha
\in(a_1^{k_1},\cdots,a_{p}^{k_p})$.
\end{lema}
\begin{prueba}
 For $j=p+1,\cdots, m$,
$\{a_1^{s_1+k_1},\cdots,a_p^{s_p+k_p},a_{p+1}^{s_{p+1}},\cdots,a_j^{s_j}\}$
is a regular $A$-sequence, and we can prove inductively that
$$\alpha a_1^{s_1}\cdots
a_{j-1}^{s_{j-1}}\in(a_1^{s_1+k_1},\cdots, a_{p}^{s_p+k_{p}}).$$
 For $i=p-1,\cdots,0$,
$\{a_1^{s_1+k_1},\cdots,a_i^{s_i+k_i},a_{i+1}^{k_{i+1}},\cdots,a_p^{k_p}\}$
is a regular $A$-sequence, and
 we inductively prove  that $$\alpha  a_1^{s_1}\cdots
a_{i}^{s_{i}}\in(a_1^{s_1+k_1},\cdots,
a_{i}^{s_i+k_i},a_{i+1}^{k_{i+1}},\cdots,a_p^{k_p}).$$
\end{prueba}

\begin{proposicion} \label{isol} Let $A$ be a commutative ring and let $I\subset A$ be
an ideal generated by a regular sequence $a_1,\dots,a_n$. Then,
the kernel of the morphism of graded algebras $$ A[X_1,\dots,X_n]
\to A[t],\quad X_i \mapsto a_it,$$ is generated by $a_iX_j
-a_jX_i$, $1\leq i< j\leq n$. In particular, $I$ is of linear
type.
\end{proposicion}

\begin{prueba}
Let $g$ be an homogeneous polynomial of degree $m$ in
$A[X_1,\cdots,X_n]$ such that $g(a_1,\cdots,a_n)=0$. Let $\exp_g =
cX^{e_g}$ be the greatest
 monomial of $g$ in the inverse lexicographic order, with
 $e_g=(s_1,\cdots,s_t,0,\cdots,0)$, $s_t\neq 0$.
  Then $$g(X_1,\cdots,X_n)-
 \exp_g \in (X_1^{s_1 +1},\cdots, X_t^{s_t +1}).$$
 By lemma \ref{regular2}, $c =\sum_{i=1}^{t-1}
\alpha_i a_i
\in
(a_1,\cdots,a_{t-1})$. Then
$$f(X_1,\cdots,X_n)=g(X_1,\cdots,X_n)-\sum_{i=1}^{t-1}
\alpha_i(a_i X_t-a_t X_i)X_1^{s_1}\cdots X_n^{s_n}$$ is an
homogeneous polynomial of degree $m$ such that $e_f<e_g$ and $$
f(X_1,\cdots,X_n)-g(X_1,\cdots,X_n) \in J=(a_i X_j-a_j X_i,
 0<i<j\leq n).$$
 In particular, $f(a_1,\cdots,a_n)=0$.
 Consequently, after a finite number of steps, we will obtain
 $h(X_1)=c_{m}X_1^m$, such
 that $h(X_1)-g(X_1,\cdots,X_n)\in J.$
  So $h(a_1)=c_{m}a_1^m=0$, $c_m=0$
and $g(X_1,\cdots,X_n)\in J.$

\end{prueba}

\section{The module $\D f^s$} \label{sub}

 Let $X$ be a $n$-dimensional complex analytic manifold, $p$ a
 point in $X$ and $f\in\hol=\hol_{X,p}$ a non zero germ of
 holomophic function with $f(p)=0$.
 Let $D$ be the  (germ of) divisor defined by $f=0$. The free module of rank one
over the ring $\hol[f^{-1},s]$ generated by the symbol $f^s$ has a
natural left module structure over the ring $\D[s]$
\cite{bern_72}: the action of a derivation $\delta\in\derhol$ is
given by $\delta (f^s) = \delta (f) s f^{-1} f^s$.

The following lemma is well-known and the proof is
straightforward.

\begin{lema} \label{coro*}
For every linear differential operator $P\in\D$ of order $d$, we
have: $$ P(f^s) = C_{P,0} f^s + C_{P,1}\binom{s}{1}f^{s-1} +
\cdots + C_{P,d}\binom{s}{d} f^{s-d}$$ where

$C_{P,d} = d!\sigma(P)(df)=\{\cdots \{\{\sigma(P),
f\},f\}\stackrel{d}{\cdots},f\}$.
\end{lema}

\medskip
Denote by $J_f\subset \hol$ the jacobian ideal associated with
$f$. The surjection $$\delta\in\derhol \mapsto \delta(f)\in J_f$$
and the canonical isomorphism of graded $\hol$-algebras
\begin{equation} \label{cano}
\simm_{\hol}(\derhol) \simeq \grD
\end{equation} induce a surjective graded morphism of $\hol$-algebras
\begin{equation} \label{varphi}
\varphi_f:\grD \longrightarrow \R(J_f).
\end{equation}

In coordinates, $\grD=\hol[\xi_1,\dots,\xi_n]$, $\xi_i =
\sigma(\partial_i)$ and $$\varphi_f(\sigma(P)) =
\sigma(P)(\partial_1(f)t, \dots,\partial_n(f)t)=\sigma(P)(df)t^d$$
for every differential operator $P\in\D$ of order $d$.

The homogeneous part of degree 1 of $\ker \varphi_f$ is naturally
identified with the $\hol$-module $$ \Theta_f =
\{\delta\in\derhol\ |\ \delta(f)=0\}$$ by means of the canonical
isomorphism (\ref{cano}).

Lemma \ref{coro*} implies that $\sigma(\Ann_{\D} f^s) \subset \ker
\varphi_f$.

\begin{proposicion} \label{nueva1}
With the above notations, if $J_f$ is of linear type, then
$$\sigma(\Ann_{\D} f^s) = \ker \varphi_f$$ and the left ideal
$\Ann_{\D} f^s$ of $\D$ is generated by $\Theta_f$.
\end{proposicion}

\begin{prueba}
By lemma \ref{kerishom}, $\ker \varphi_f = \grD \Theta_f
\subset\sigma(\Ann_{\scr D}f^s)$.

\smallskip The inclusion $\D \Theta_f \subset \Ann_{\D} f^s$ is
obvious. Let's prove that $\Ann_{\D} f^s \subset \D \Theta_f$.
Clearly, $F^1 \Ann_{\D} f^s =\Theta_f$. Suppose
$F^{d-1}\Ann_{\D}f^s \subset \D \Theta_f$ and take a differential
operator $P \in F^d\Ann_{\scr D}f^s \setminus F^{d-1}\Ann_{\D}f^s
$. Then, $\sigma(P)\in \ker \varphi_f=\grD \sigma(\Theta_f)$, and
$\sigma(P)=\sum A_i\sigma(\delta_i),$ where $\delta_i\in \Theta_f$
and the $A_i$ are homogeneous of degree $d-1$. Let $Q_i$ be
differential operators such that $\sigma(Q_i)=A_i$. We apply the
induction hypothesis to $P - \sum_iQ^i\delta_i\in
F^{d-1}\Ann_{\scr D}f^s$ and we conclude the result.
\end{prueba}

\begin{proposicion} \label{isosin}
 {\em (Isolated singularities case, cf. \cite{mai_effec}, 2.7)}
 If $f$ has isolated singularity, then $\ker \varphi_f$ is
generated by $\partial_i(f)\xi_j -\partial_j(f)\xi_i$, $1\leq i<
j\leq n$. In particular, the left ideal $\Ann_{\D} f^s$ of $\D$ is
generated by $\partial_i(f)\partial_j -\partial_j(f)\partial_i$,
$1\leq i< j\leq n$.
\end{proposicion}

\begin{prueba} It is a consequence of lemma \ref{isol} and proposition \ref{nueva1}.
\end{prueba}

\bigskip
\section{Locally quasi-homogeneous free divisors are \\ Koszul free}

\begin{proposicion} \label{uf}
Let $U$ be a connected open set of a complex $n$-dimensional
analytic manifold $X$  and let $\Sigma\subset U$ be a closed
analytic set of dimension $s$. If a sequence
$C=\{\sigma_1,\dots,\sigma_{n-s}\}$ of homogeneous polynomials in
$\OX(U)[\xi_1,\dots,\xi_n]$ is regular at every point $q\in
U\setminus \Sigma$ (i.e. it is regular in
$\hol_{X,q}[\xi_1,\dots,\xi_n]$), then it is regular at every
point of $U$.
\end{proposicion}

\begin{prueba}
Let $p\in \Sigma$ and let $\pi:U\times\CC^n\to U$ be the
projection. By proposition \ref{equiv}, we have to prove that the
ideal $I=(\sigma_1,\cdots,\sigma_{n-s})$ defines an analytic set $
V=V(I) \subset U\times\CC^n$ of dimension $n+s$. By hypothesis, we
know that $C$ is regular on $U\setminus \Sigma$, and so ({\em
loc.~cit.}) the dimension of (every irreducible component of)
$V\cap\pi^{\-1}(U\setminus\Sigma)$ is $n+s$. Now, let $W$ be an
irreducible component of $V$. It has, at least, dimension $n+s$.
 If $W$ is contained in $\pi^{-1}(\Sigma)=\Sigma\times\CC^n$, then it must be equal to
 $\pi^{-1}(\Sigma)$.
 If not, $\dim W = \dim (W\cap \pi^{-1}(U\setminus\Sigma)) \leq
\dim
 (V\cap\pi^{-1}(U\setminus\Sigma)) = n+s$. So, we conclude that
$W$ has dimension $n+s$.
\end{prueba}

\begin{corolario} \label{uff}
Let $D$ be a free divisor in some analytic manifold $X$ and let
$\Sigma\subset D$ a discrete set of points. If $D$ is Koszul free
at every point $x\in D\setminus \Sigma$, then $D$ is Koszul free
(at every point of $X$).
\end{corolario}

\begin{teorema} \label{lqhiskfd} Every locally quasi-homogeneous
free divisor is Koszul free.
\end{teorema}
\begin{prueba}
We proceed by induction on the dimension $t$ of the ambient
manifold $X$. For $t=1$, the theorem is trivial and for $t=2$, the
theorem is directly proved in example \ref{examples}, 3).
 Now, we suppose that the result is
true for $t < n$, and let $D$ be a locally quasi-homogeneous free
divisor of a complex analytic manifold $X$ of dimension $n$.
 Let $p\in D$ and let  $\{\delta_1,\dots,\delta_n
 \}$ be a basis of the logarithmic derivations of $D$ at $p$.

Thanks to \cite{cas_mon_nar_96}, prop. 2.4 and lemma 2.2,~(iv),
there is an open neighborhood $U$ of $p$ such that for each $q\in
U\cap D$, with $q\neq p$, the germ of pair $(X,D,q)$ is isomorphic
to a product $ (\CC^{n-1}\times \CC, D^{\prime}\times \CC,
(0,0))$, where $D^{\prime}$ is a locally quasi-homogeneous free
divisor. Induction hypothesis implies that $D^{\prime}$ is a
Koszul free divisor at $0$. Then, by proposition \ref{bases}.1.,
$D$ is a Koszul free divisor at $q$ too. We have then proved that
$D$ is a Koszul free divisor in $U\setminus \{p\}$. We conclude by
using corollary \ref{uff}.
\end{prueba}

\begin{corolario} Every free divisor that is locally
quasi-homogeneous at the complement of a discrete set is Koszul
free.
\end{corolario}

In particular, the last corollary gives rise a new proof of the
fact that every divisor in dimension 2 is Koszul free (cf.
\ref{examples}, 3)).

\section{The module ${\scr D}f^s$ for locally quasi-homogeneous
free divisors}

\numero \label{entorno}  In this section, $f\in\hol = \hol_{X,p}$
will be a reduced locally quasi-homogeneous free germ \ref{def-1},
\ref{kfd}. That means that $D=\{f=0\}$ is a locally
quasi-homogeneous free divisor near $p$.

We will also assume that
\begin{enumerate}
\item[-)] The equation
$f$ and its Euler vector field $E$ are globally defined on $X$.
\item[-)] $E(q)\neq 0$ for every $q\in X\setminus\{ p\}$.
\item[-)] $\derD$ is $\OX$-free (of rank $n=\dim X$).
\end{enumerate}

\medskip
In order to proceed inductively on the dimension of the ambient
variety when working with such $f$'s, we quote the following
direct consequence of \cite{MR80i:32027}, lemmas 1.3, 1.5 (see
also \cite{cas_mon_nar_96}, prop. 2.4)

\begin{proposicion} \label{preparacion-2} Let $f\in \hol_{X,p}$ a reduced locally
quasi-homogeneous free germ and let $D$ be the divisor $f=0$. For
$q\in D\setminus \{p\}$ close to $p$, there are local coordinates
$z_1,\dots,z_n\in\hol_{X,q}$ centered at $q$ and a
quasi-homogeneous polynomial $G'(t_1,\dots,t_{n-1})$ in $n-1$
variables which is also a locally quasi-homogeneous free germ in
$\hol_{\CC^{n-1},0}$ and  such that $f_q = G'(z_1,\dots,z_{n-1})$.
\end{proposicion}

We call $\widetilde{\Theta}_f$ the $\OX$-sub-module (and Lie
algebra) of $\derD$ whose sections are vector fields annhilating
$f$. Denote by $\J_f\subset \OX$ the jacobian ideal sheaf
associated with $f$. The stalk of $\widetilde{\Theta}_f$ (resp. of
$\J_f$) at $p$ is then $\Theta_f$ (resp. $J_f$).

As in (\ref{varphi}), we have a surjective graded morphism of
$\OX$-algebras $$\Phi_f:\GrD \longrightarrow \R(\J_f),$$ whose
stalk at $p$ is $\varphi_f$.

We have:
\begin{equation}\label{sumadir}
\derD = \widetilde{\Theta}_f \oplus (\OX E),\quad \derlogf =
\Theta_f \oplus (\hol E),
\end{equation}
and $\widetilde{\Theta}_f, \Theta_f$ are free of rank $n-1$.

\begin{proposicion} \label{KC}
The Koszul complex associated with $\widetilde{\Theta}_f \subset
\DER = \Gr^1_{F^{\bullet}}(\DX) \subset \GrD $:
\begin{eqnarray*}
0 \xrightarrow{}
\GrD\otimes_{\OX}\stackrel{n-1}{\bigwedge}\widetilde{\Theta}_f
\xrightarrow{d_{-n+1}} \cdots \xrightarrow{d_{-2}}
\GrD\otimes_{\OX}\stackrel{1}{\bigwedge}\widetilde{\Theta}_f
\xrightarrow{d_{-1}} \GrD \\ d_{{-k}}(
F\otimes(\sigma_1\wedge\cdots\wedge\sigma_k))
 =\sum_{i=1}^k
(-1)^{i-1} P\sigma_i\otimes(\sigma_1\wedge\cdots\widehat{\sigma_i}
\cdots\wedge\sigma_k), \quad 0 < k < n,
\end{eqnarray*}
is exact.
\end{proposicion}

\begin{prueba} We need to prove that some (or any) basis
$\{\delta_1,\dots,\delta_{n-1}\}$ of $\widetilde{\Theta}_f$ form a
regular sequence in $\GrD$, but such a basis can be augmented to a
basis $\{\delta_1,\dots,\delta_{n-1},E\}$ of $\derD$, that we know
by theorem \ref{lqhiskfd} to form a regular sequence in $\GrD$.
\end{prueba}

\begin{proposicion} \label{localcohom}
With the hypothesis of \ref{entorno}, if the augmented graded
complex of $\GrD$-modules
\begin{equation} \label{CKF}
0\to\GrD\otimes_{\OX}\stackrel{n-1}{\bigwedge}\widetilde{\Theta}_f
\xrightarrow{d_{-n+1}} \cdots \xrightarrow{d_{-1}} \GrD
\xrightarrow{\Phi_f} \R(\J_f)\xrightarrow{} 0
\end{equation}
 is exact on $X-\{p\}$, then it is
exact everywhere.
\end{proposicion}

\begin{prueba} We know that $\Phi_f$ is surjective. By proposition \ref{KC}, the only thing to
prove is $\ker \Phi_ f = \IM d_{-1}$. We can proceed separately on
each homogeneous component:
$$0\to \Gr^{m-n+1}_{F^{\bullet}}(\DX)
\otimes_{\OX}\stackrel{n-1}{\bigwedge}\widetilde{\Theta}_f
\xrightarrow{d_{-n+1}^m} \cdots \xrightarrow{d_{-1}^m}
\Gr^m_{F^{\bullet}}(\DX) \xrightarrow{\Phi_f^m} \J_f^m
\xrightarrow{} 0.$$ Let's consider the coherent $\OX$-module
${\scr F}=\IM d_{-1}^m$ and the short sequence
\begin{equation} \label{ss}
 0 \xrightarrow{} {\scr F} \xrightarrow{}
\Gr^m_{F^{\bullet}}(\DX) \xrightarrow{\Phi_f^m} \J_f^m
\xrightarrow{} 0.
\end{equation}
By proposition  \ref{KC} and the fact that the cohomology with
support $H^i_p({\OX})$ vanishes for $i\neq n$, we deduce that
$H^i_p({\scr F})=0$ for $i=0,1$ and $ H^0_p({\J_f^m})=0$. These
properties and the exactness of (\ref{ss}) on $X-\{p\}$ imply the
proposition (cf. \cite{loo_84}, (8.14)).
\end{prueba}

 The following lemma is clear.

\begin{lema} \label{DyDprima}
Let $g\in \hol_{n-1}=\CC\{y_1,\cdots,y_{n-1}\}$ and call $f=g$,
but as an element in $\hol_n=\CC\{y_1,\cdots,y_n\}$. Then:
\begin{enumerate}
\item $\ker \varphi_f$ is generated by $\ker\varphi_g$ and
$\sigma({\partial_{y_n}})$.
\item $\Theta_f$ is generated by $\Theta_g$ and ${\partial_{y_n}}$.
\end{enumerate}
\end{lema}

\begin{teorema} \label{teorema}
Let $f\in\hol = \hol_{X,p}$ be a reduced locally quasi-homogeneous
free germ. Then the graded complex of $\grD$-modules
$$0\to\grD\otimes_{\hol}\stackrel{n-1}{\bigwedge}{\Theta_f}
\stackrel{\varepsilon_{-n+1}} {\to}\cdots
\stackrel{\varepsilon_{-1}} {\to} \grD \stackrel{\varphi_f}{\to}
\R(J_f)\to 0$$ is exact. In particular, the kernel of the morphism
$$\grD \stackrel{\varphi_f}{\to} \R(J_f)$$ is the ideal generated
by $\Theta_f$ and then the jacobian ideal $J_f$ is of linear type.
\end{teorema}

\begin{prueba}
By the exactness of (\ref{KC}), the only thing to prove is that
$\ker \varphi_f$ is generated by $\sigma(\Theta_f)$. We will use
induction on $n=\dim X$. If $n=2$, we apply proposition
\ref{isosin}.
 We suppose that the result is true if the ambient
variety has dimension $n-1$.
 By Proposition \ref{localcohom},
we need to prove the exactness of the complex (\ref{CKF}) on
$U\setminus\{x\}$, for some open neighborhood $U$ of $x$, or
equivalently, that $\ker \Phi_f$ is generated by
$\sigma(\Theta_f)$ at every $q\in U\setminus \{x\}$.
 The result is then a consequence of proposition
 \ref{preparacion-2}, lemma \ref{DyDprima} and  the induction hypothesis.
\end{prueba}

\begin{definicion} \label{spen-2}
The Spencer complex\footnote{It should be noticed that such
complex was originally used by Chevalley and Eilenberg in the
setting of the cohomology of Lie algebras (cf. \cite{wei_94},
7.7).}
 for $\widetilde{\Theta}_f$ is the
complex of free left $\D_X$-modules given by: $$ 0 \xrightarrow{}
\DX\otimes_{\OX}\stackrel{n-1}{\bigwedge}\widetilde{\Theta}_f
\xrightarrow{\varepsilon_{-n+1}} \cdots
\xrightarrow{\varepsilon_{-2}}
\DX\otimes_{\OX}\stackrel{1}{\bigwedge}\widetilde{\Theta}_f
\xrightarrow{\varepsilon_{-1}} \DX,$$ $$\varepsilon_{-1} (
 P\otimes\delta)=P\delta; \quad
\varepsilon_{{-k}}( P\otimes(\delta_1\wedge\cdots\wedge\delta_k))
 =\sum_{i=1}^k
(-1)^{i-1} P\delta_i\otimes(\delta_1\wedge\cdots\widehat{\delta_i}
\cdots\wedge\delta_k)$$ $$+ \sum_{1\leq i<j\leq k}
(-1)^{i+j}P\otimes([\delta_i,\delta_j]\wedge\delta_1\wedge\cdots
\widehat{\delta_i}\cdots\widehat{\delta_j}
       \cdots\wedge\delta_k), \ \ 2\leq k\leq n-1. $$
In a similar way we define the Spencer complex for $\Theta_f$,
which is the stalk at $p$ of the Spencer complex for
$\widetilde{\Theta}_f$.

Both Spencer complexes can be augmented by considering the obvious
maps $\DX \to \DX f^s, \D \to\D f^s$.
\end{definicion}

\begin{corolario}  \label{coro-x}
With the hypothesis of \ref{entorno}, we have:
\begin{enumerate}
\item[(a)]
 The Spencer complex
for ${\Theta}_f$ is a resolution of  $\D f^s$. In particular, the
left ideal $\Ann_{\D} f^s$ is generated by ${\Theta}_f$.
\item[(b)] The left ideal $\Ann_{\D[s]} f^s$ is generated by
${\Theta}_f$ and $E-s$.
\item[(c)] The left ideal $\Ann_{\D} \eta$, where $\eta$ is the class of
$f^s$ in the quotient $\D[s]f^s/\D[s]f^{s+1}$, is generated by
${\Theta}_f$ and $f$.
\end{enumerate}
\end{corolario}

\begin{prueba} For (a) we proceed as in \cite{cal_99}, prop. 4.1.3 by
using proposition \ref{nueva1} and theorem \ref{teorema}. Property
(b) follows easily from (a), and property (c) follows from (a) and
(b).
\end{prueba}

 Let's call ${\Xi}_f = {\Theta}_f \oplus (\hol f)$ (resp. $\widetilde{\Xi}_f =
\widetilde{\Theta}_f \oplus (\OX f)$), which is a free
sub-$\hol$-module (respectively, sub-$\OX-$module) and a Lie
subalgebra of $\D$ (resp. of $\DX$). It can be also canonically
embedded in $\grD$ (resp. $\GrD$) equipped with the Poisson
bracket $\{-,-\}$. As in \ref{KC} and \ref{spen-2}, we define the
Koszul complex associated with
 ${\Xi}_f\subset \grD$ (resp. $\widetilde{\Xi}_f\subset \GrD$) and the
Spencer complex associated with ${\Xi}_f \subset \D$ (resp.
$\widetilde{\Xi}_f \subset \DX$). The Koszul (resp. Spencer)
complex associated with ${\Xi}_f\subset \grD$ (resp. with ${\Xi}_f
\subset \D$) is obviously the stalk at $p$ of the Koszul (resp. of
the Spencer) complex associated with $\widetilde{\Xi}_f\subset
\GrD$ (resp. with $\widetilde{\Xi}_f \subset \DX$)

\begin{teorema} \label{mainn}
With the hypothesis of \ref{entorno}, the following properties hold:
\begin{enumerate}
\item The Koszul complex associated with $\Xi_f\subset \grD$ is
exact.
\item The Spencer complex associated with $\Xi_f \subset \D$ is
a free resolution of $\D[s]f^s/\D[s]f^{s+1}$.
\end{enumerate}
\end{teorema}

\begin{prueba} For the first property, call $\bK$ the Koszul complex
associated with $\widetilde{\Xi}_f\subset \GrD$. The  Koszul
complex associated with $\Xi_f\subset \grD$ is the stalk at $p$ of
$\bK$.

We proceed by induction on the dimension of the ambient variety.
If that dimension is 1, $\Xi_f = \hol f$, and the Koszul complex
associated with $f$ is clearly exact. Suppose the result true if
the dimension of the ambient variety is $< n$.

Now, suppose $\dim X = n$.

Let $\delta_1,\dots,\delta_{n-1}$ be a basis of
$\widetilde{\Theta}_f$ in some small enough neighborhood $U$ of
$p$. According to proposition \ref{uf}, we need to prove that
$\bK$ is exact on $U\setminus \{p\}$.

For every $q\in U$ with $f(q)\neq 0$, the germ of $f$ at $q$ is an
unit and by proposition \ref{KC}, the complex $\bK$ is exact at
$q$.

Let $q$ be a point in $D = \{f=0\}$, $q\neq p$. By proposition
\ref{preparacion-2}, there are local coordinates
$z_1,\dots,z_n\in\hol_{X,q}$ and a quasi-homogeneous polynomial
$G'(t_1,\dots,t_{n-1})\in\hol_{\CC^{n-1},0}$ in $n-1$ variables
which is also a locally quasi-homogeneous free germ in
$\hol_{\CC^{n-1},0}$ and such that $f_q = G'(z_1,\dots,z_{n})$.

Let $G(t_1,\dots,t_n)\in\hol_{\CC^{n},0}$ be the same polynomial
as $G'(t_1,\dots,t_{n-1})$ but considered  in $n$ variables. The
exactness of $\bK_q$ is then equivalent to the exactness of the
Koszul complex associated with $\Xi_G\subset\gr_{F^{\bullet}}
\D_{\CC^{n},0}$.

Let's write $\hol_m = \CC\{t_1,\dots,t_m\}$ and call $\xi'_i$ the
principal symbol of $\frac{\partial}{\partial t_i}$.

Let $\{\delta'_1,\dots,\delta'_{n-2}\}\subset
\oplus_{i=1}^{n-1}\hol_{n-1}\frac{\partial}{\partial t_i}$ be a
basis of $\Theta_{G'}$. A basis of $\Theta_G$ is then
$\{\delta'_1,\dots,\delta'_{n-2},\frac{\partial}{\partial
t_n}\}\subset \oplus_{i=1}^{n}\hol_{n}\frac{\partial}{\partial
t_i}$. Call  $\sigma'_i$ the principal symbol of $\delta'_i,
i=1,\dots,n-2$.

By induction hypothesis we know that the Koszul complex associated
with $\Xi_{G'} \subset \gr_{F^{\bullet}}
\D_{\CC^{n-1},0}=\hol_{n-1}[\xi'_1,\dots,\xi'_{n-1}]$ is exact or,
equivalently, that $\sigma'_1,\dots,\sigma'_{n-2},G'$ is a regular
sequence in $\hol_{n-1}[\xi'_1,\dots,\xi'_{n-1}]$. That implies
that $\sigma'_1,\dots,\sigma'_{n-2},\xi'_n,G=G'$  is a regular
sequence in $\hol_n[\xi'_1,\dots,\xi'_n]$, i.e. that the Koszul
complex associated with $\Xi_{G} \subset \gr_{F^{\bullet}}
 \D_{\CC^{n},0}$ is
exact, and the result is proved.

\smallskip
For the second property, we filter the Spencer complex associated
with $\Xi_f \subset \D$ as in \cite{gros_nar}, prop.~2.3.4: $$
\deg (\Theta_f) = 1, \deg (f) = 0.$$ Its graded complex coincides
with the Koszul complex associated with $\Xi_f\subset \grD$, and
then the Spencer complex is exact. To conclude, we use corollary
\ref{coro-x}, (c).
\end{prueba}

\section{Examples and questions}

 We know several (related) kind of free divisors:
 \begin{enumerate}
 \item[{[LQH]}] Locally quasi-homogeneous (definition \ref{lqh}).
 \item[{[EH]}] Euler homogeneous (definition \ref{eh}).
 \item[{[LCT]}] Free divisors satisfying the logarithmic comparison
 theorem.
 \item[{[KF]}] Koszul free (definition \ref{kfd}).
 \item[{[P]}] Free divisors such that the complex $\Olog$ is
 a perverse sheaf.
 \end{enumerate}

We have then the following implications: \begin{center} [LQH]
$\Rightarrow$ [EH] (obvious), \hspace{0.5cm}
 {[LQH]} $\Rightarrow$ {[LCT]}  by \cite{cas_mon_nar_96}, th. 1.1,\\
 {[LCT]} $\Rightarrow$ [P], by \cite{mebkhout}, II, th. 2.2.4)  \hspace{0.5cm}
 {[KF]} $\Rightarrow$ [P] by \cite{cal_99}, th. 4.2.1, \\
 {[LQH]} $\Rightarrow$ {[KF]} by theorem \ref{lqhiskfd}.\end{center}

\begin{ejemplo} (Free divisors in dimension 2)
We recall theorem 3.9  from \cite{cal_mon_nar_cas}: Let $X$ be a
complex analytic manifold of dimension 2 and $D\subset X$ a
divisor. The following conditions are equivalent:\\ 1. $D$ is
Euler homogeneous.\\ 2. $D$ is locally quasi-homogeneous.\\ 3. The
logarithmic comparison theorem holds for $D$.

Consequently, in dimension 2 we have: $$ \text{[LQH]}
\Leftrightarrow \text{[EH]} \Leftrightarrow \text{[LCT]}$$ and
[KF] (cf. \ref{examples}, 3) and [P] (\cite{cal_99}) always hold.
In particular, $$\text{[KF]} \not\Rightarrow \text{[LQH], [EH],
[LCT]}.$$

Examples of plane curves not satisfying logarithmic comparison
theorem are, for instance, the curves of the family (cf.
\cite{cal_mon_nar_cas}): $$x^q+y^q+xy^{p-1}=0,\quad p\geq q+1\geq
5.$$
\end{ejemplo}

\begin{ejemplo} \label{rectas}
(An example in dimension 3) Let's consider $X=\CC^3$ and $D
=\{f=0\}$, with $f=xy(x+y)(x+yz)$ \cite{cal_99}. A basis of
$\derD$ is $\{\delta_1,\delta_2,\delta_3\}$, with $$
\begin{array}{ccrcrcr}
\delta_1 &=& x y \Dx &+& y^2 \Dy &-&4 (x + y z) \Dz,\\ \delta_2
&=&  x(x+3y)
  \Dx &-& y (3x+y) \Dy &+& 4 x (z-1) \Dz,\\
\delta_3&=&x\partial_x &+&y \partial_y
\end{array}
$$ the determinant of the coefficients matrix being $-16f$ and
$$\delta_1(f) = 0,\quad\delta_2(f) = 0,\quad \delta_3(f) = 4f.$$

In particular, $D$ is Euler homogeneous ($E = (1/4) \delta_3$) and
we know \cite{cal_mon_nar_cas} that it satisfies the logarithmic
comparison theorem.
 Let $I\subset {\scr O}_{T^{\ast}X}$ be the ideal generated by the
 symbols $\{\sigma_1,\sigma_2,\sigma_3\}$ of the basis of $\derD$.
By corollary \ref{2.9}, $D$ is not Koszul free, because the
dimension of $V(I)$ at $((0,0,\lambda), 0)\in T^{\ast}X$ is 4, and
neither is $D$ locally quasi homogeneous.

\medskip
So:
\begin{center}
[LCT]  $\not\Rightarrow$ [KF], [LQH],\hspace{0.5cm} [EH]
$\not\Rightarrow$ [KF], [LQH].\end{center}
\end{ejemplo}

Finally, for the only missing relation,
we quote the following
conjecture from \cite{cal_mon_nar_cas}:

 \begin{conjetura}
 If the logarithmic comparison theorem holds for $D$, then $D$ is
 Euler homogeneous.
\end{conjetura}

\begin{ejemplo} Let's see in the
 example \ref{rectas} that  the left ideal $\Ann_{\D}(f^s)$ is not
 generated by $\Theta_f$ and then, $J_f$ is not an ideal
 of linear type.

\smallskip Here, we set $X=\CC^3$, $p = (0,0,0)$ and $E=(1/4)\delta_3$.
  The $\hol$-modules  $\Theta_f$ and
  $\derlogf= \Theta_f \oplus \hol \cdot E $
 are generated  by $\{\delta_1, \delta_2\}$
 and $\{\delta_1,\delta_2,E\}$, respectively.
 The  symbols $\sigma_1=\sigma(\delta_1)$,
 $\sigma_2=\sigma(\delta_2)$ form a $\grD$-regular sequence (the proof is analogous to
  example \ref{examples}, 3)). Then, as in the proof of
  \cite{cal_99}, prop. 4.1.2, we have
  $$\sigma (\D\Theta_f) = \grD \sigma(\Theta_f).$$

 For
\begin{eqnarray*}
P=2y^2\Dx\Dy-2y^2\Dy^2-   (2xz+6yz)\Dx\Dz+10yz\Dy\Dz +8z(1 -z)\Dz^2+\\
+ (2x-4y)\Dy\Dz-x\Dx-y\Dy-8z\Dz+4\Dz
 \end{eqnarray*}
 and $R =\CC[x,y,z]$, $S = R[\xi_1,\xi_2,\xi_3]$, ${\mathfrak m} =
 R(x,y,z)$
 we check that
 \begin{enumerate}
\item $P\in  \Ann_{\DX}(f^s)$,
\item $(S(\sigma_1,\sigma_2):\sigma(P)) = S(x,y)$, and then
$(S(\sigma_1,\sigma_2):\sigma(P))\cap R = R(x,y)$.
\end{enumerate}

 So, $\sigma(P)\notin R_{\mathfrak m}[\xi_1,\xi_2,\xi_3]
\sigma(\Theta_f)$ and, by faithful flatness, $$\sigma (P) \notin
\hol [\xi_1,\xi_2,\xi_3] \sigma(\Theta_f) = \grD
\sigma(\Theta_f).$$ We conclude that $ P \notin \D\Theta_f.$
\end{ejemplo}

\begin{problema} We do not know whether a free divisor defined by
a quasi-homogeneous polynomial (with strictly positive weights) is
locally quasi-homogeneous.
\end{problema}

\begin{problema} We do not know any example of a free divisor
$D\subset X$ whose logarithmic de Rham complex $\Olog$ is not
perverse.
\end{problema}

{\small
\bibliographystyle{plain}

\begin{thebibliography}{10}

\bibitem{alek_90}
A.G. Aleksandrov.
\newblock Nonisolated hypersurface singularities.
\newblock In {\em Theory of singularities and its applications}, volume~1 of
  {\em Adv. Soviet Math.}, pages 211--246. A.M.S., Providence, R.I., 1990.

\bibitem{bern_72}
J.~Bernstein.
\newblock The analytic continuation of generalized functions with respect to a
  parameter.
\newblock {\em Funz. Anal. Appl.}, 6 (1972), 26--40.

\bibitem{calde_tesis}
F.J.~Calder{\'o}n~Moreno.
\newblock Operadores diferenciales logar\'{\i}tmicos con respecto a un divisor
  libre.
\newblock Univ. Sevilla, July 1997.
\newblock Ph.D.

\bibitem{cal_99}
F.J.~Calder\'on-Moreno.
\newblock Logarithmic Differential Operators and Logarithmic De Rham
    Complexes Relative to a Free Divisor.
\newblock {\it Ann. Sci. E.N.S.}, 32 (1999), 577-595.

\bibitem{cal_mon_nar_cas}
F.J.~Calder\'{o}n-Moreno, D.Q. Mond, L.~Narv\'{a}ez-Macarro and
F.J.~Castro-Jim\'{e}nez.
\newblock Logarithmic Cohomology of the Complement of a Plane Curve.
\newblock {\it  Comment. Math. Helv.}, 77 (2002), 24-38.

\bibitem{lqhf-kf}
F.J. Calder{\'o}n-Moreno and L. Narv{\'a}ez-Macarro.
\newblock Locally quasi-homogeneous free divisors are {K}oszul free.
\newblock {\it Proc. Steklov Inst. Math.}, 238 (2002), to appear.

\bibitem{cas_mon_nar_96}
F.J.~Castro-Jim\'{e}nez, D.~Mond and L.~Narv\'{a}ez-Macarro.
\newblock Cohomology of the complement of a free divisor.
\newblock {\it Transactions of the A.M.S.}, 348 (1996), 3037--3049.

\bibitem{deligne} P.~Deligne. \'{E}quations Diff\'{e}rentielles \`{a}
Points Singuliers R\'{e}guliers, {\it Lect. notes in Math., 163},
Springer-Verlag, Berlin-Heidelberg, 1970.

\bibitem{MR80i:32027}
R.~Ephraim.
\newblock Isosingular loci and the {C}artesian product structure of complex
  analytic singularities.
\newblock {\em Trans. Amer. Math. Soc.}, 241 (1978), 357--371.

\bibitem{gros_nar}
M.~Gros and L.~Narv\'aez-Macarro.
\newblock Cohomologie {\'e}vanescente p-adique: calculs locaux.
\newblock {\em Rendiconti Sem. Mat. Univ. Padova}, 104 (2000), 71--90.

\bibitem{grot} A. Grothendieck. \newblock On the de Rham cohomology of
algebraic varieties. \newblock {\it Publ. Math. de l'I.H.E.S.}, 29
(1966), 95-103.

\bibitem{loo_84}
E.J.N.~Looijenga.
\newblock {\em Isolated singular points on complete intersections},
 {\em London Mathem. Soc. Lect. Notes Series, 77}.
\newblock Cambridge Univ. Press, Cambridge, 1984.

\bibitem{mai_effec}
Ph.~Maisonobe.
\newblock {$\scr D$}-modules: an overview towards effectivity.
\newblock In {\em Computer algebra and differential equations}, volume 193 of
  {\em London Math. Soc. Lecture Note Ser.}, pages 21--55. Cambridge Univ.
  Press, 1994.

\bibitem{mat} H.~Matsumura. Commutative ring theory, Cambridge
University Press, Cambridge, 1992.

\bibitem{mebkhout} Z.~Mebkhout. Le formalisme des six op\'{e}rations de
Grothendieck
 pour les $\DX$-modules coh\'{e}rents,
 {\it Travaux en cours, 35}, Hermann, Paris, 1989.

\bibitem{saito}
{K.~Saito.}
\newblock On the uniformization of complements of discriminant
loci.
\newblock Preprint, Williams College, 1975.

\bibitem{ksaito_log}
K.~Saito.
\newblock Theory of logarithmic differential forms and logarithmic vector
  fields.
\newblock {\em J. Fac. Sci. Univ. Tokyo}, 27 (1980), 265--291.

\bibitem{vascon_cmcaag}
W.V.~Vasconcelos.
\newblock {\em Computational methods in commutative algebra and algebraic
  geometry}, {\em Algorithms and Computation in Mathematics, 2}.
\newblock Springer Verlag, New York, 1998.

\bibitem{wei_94}
Charles~A. Weibel.
\newblock {\em An introduction to homological algebra}, volume~38 of {\em
  Cambridge studies in advanced mathematics}.
\newblock Cambridge University Press, Cambridge, 1994.

\end{thebibliography}

}

{\small Departamento de \'{A}lgebra,
 Facultad de  Matem\'{a}ticas, Universidad de Sevilla, P.O. 1160, 41080
 Sevilla, Spain}. \\
{\small {\it E-mail}: \begin{itemize} \item[]
calderon@algebra.us.es
 \item[] narvaez@algebra.us.es
\end{itemize} }

\end{document}